\documentstyle[11pt,amssymb,amsfonts]{article}
  
\begin{titlepage}

\title{On the functoriality of  Lott's secondary analytic index}

\author{Ulrich Bunke\thanks{Mathematisches Institut, Universit\"at G\"ottingen, Bunsenstr. 3-5, 37073 G\"ottingen, GERMANY, 
E-mail:bunke@uni-math.gwdg.de}}

\end{titlepage}

\newcommand{\C}{{\bf C}}

\newcommand{\image}{{\rm image}}

\newcommand{\Tr}{{\rm Tr}}
\newcommand{\cT}{{\cal T}}

\newcommand{\cH}{{\cal H}}

\newcommand{\cE}{{\cal E}}
\newcommand{\ccR}{{\cal R}}
\newcommand{\cV}{{\cal V}}

\newcommand{\cC}{{\cal C}}

\newcommand{\Hom}{{\mbox{\rm Hom}}}

\newcommand{\End}{{\mbox{\rm End}}}

\newcommand{\cF}{{\cal F}}
\newcommand{\cD}{{\cal D}}

\def\hB{\hspace*{\fill}$\Box$ \newline\noindent}

\newcommand{\cS}{{S}}

\def\hB{\hspace*{\fill}$\Box$ \\[0.5cm]\noindent}

\newtheorem{prop}{Proposition}[section]

\newtheorem{theorem}[prop]{Theorem}

\begin{document}
\maketitle
\tableofcontents

\section{Introduction}

In \cite{lott99} for a smooth manifold $M$ J. Lott has defined a secondary
$K$-group $\bar{K}_R^0(M)$ of local systems of $R$-modules ($R$ is some ring)
such that their complexifications (we must fix a representation
$\rho:R\rightarrow \End(\C^n)$ in order to define  this notion) have
explicitly trivialized characteristic classes. If $p:E\rightarrow B$ is a
smooth fibre bundle with closed fibres, then he defines a push-forward
$p_!:\bar{K}_R^0(E)\rightarrow \bar{K}_R^0(B)$. This construction involves the
analytic torsion form of \cite{bismutlott95}. We will review Lott's
construction in Section \ref{lott}.   In \cite{ma99}  X. Ma has studied the
behaviour of the analytic torsion form under iterated fibre bundles. We review
his result in Section \ref{ma}.

The goal of the present note is to verify that Ma's result indeed implies
that the push-forward of Lott is functorial. We now describe this assertion in
detail. Let $p:E\rightarrow B$ be a fibration with fibre $Z$ which is in fact
an iterated fibration. So we assume that there are fibrations 
$p_1:E=:E_1\rightarrow E_2$ with fibre $Z_1$ and
$p_2:E_2\rightarrow B$ with fibre $Z_2$ such that $Z\rightarrow Z_2$ is a
fibration with fibre $Z_1$. We will verify the following:
\begin{theorem}\label{main}
We have
$p_!=(p_2)_!\circ (p_1)_!$ as maps from $\bar{K}_R^0(E)$ to $\bar{K}_R^0(B)$.
\end{theorem} 
The details of the proof are contained in Section \ref{pf}.  
{\em We do now claim any kind of originality here because the argument just
amounts to combine the construction of Lott and the result of Ma
mentioned above.}

\section{Secondary $K$-Groups and analytic push-forward}\label{lott}

This section reviews Sections 2.1 and 2.2 of \cite{lott99}. 
Suppose that $R$ is a right-Noetherian ring which is right-regular,
i.e. every finitely generated right-$R$-module has a finite resolution by
finitely generated projective right-$R$-modules. We fix a representation
$\rho:R\rightarrow \End(\C^n)$ such that $\C^n$ becomes a flat left $R$-module.
If $V$ is a right-$R$-module, then $V_\C:=V\otimes_{R,\rho} \C^n$ is called its complexification.

Let $M$ be a connected smooth manifold.
If $\cF$ is a local system (a locally constant sheaf) of finitely generated right-$R$-modules, 
then $\cF_\C$ (the fibrewise complexification) is the sheaf of parallel
sections of a flat complex vector bundle which we denote
by $(F_\C,\nabla^{F_\C})$.

Let $(E,\nabla^E)$ be a flat complex vector bundle. 
If we choose a hermitean metric $h^E$ (we consider $h^E$ is a section of the
flat bundle $\Hom_\C(E,E^*)$)  on $E$, then we can define the characteristic
forms 
\begin{eqnarray*}
\omega(\nabla^{E},h^{E})&:=&(h^E)^{-1} \nabla^{\Hom(E,E^*)}h^{E}\\
c_k(\nabla^{E},h^{E})&=&(2\pi\imath)^{-\frac{k-1}{2}} 2^{-k} \Tr\:
\omega(\nabla^{E},h^{E})\\
c(\nabla^{E},h^{E})&=&\sum_{j=0}^\infty
\frac{1}{j! }c_{2j+1}(\nabla^{E},h^{E})\ . \end{eqnarray*}
The form $c(\nabla^{E},h^{E})$ is closed and represents the characteristic
class $c(\nabla^E)\in H_{dR}^{odd}(M)$ of the flat vector bundle $
(E,\nabla^E)$ which is independent of the choice of $h^E$.

The abelian group $\hat{K}^0_R(M)$ is generated by triples
$f=(\cF,h^{F_\C},\eta)$, where 
\begin{enumerate}
\item $\cF$ is a local system of finitely generated right-$R$-modules,
\item $h^{F_\C}$ is a hermitean metric of the corresponding flat complex
vector bundle  $(F_\C,\nabla^{F_\C})$, and 
\item $\eta\in\Omega^{ev}(M)/\image(d)$,
\end{enumerate}
subject to the following relations :
If 
$$0\rightarrow \cF_1\rightarrow \cF_2\rightarrow \cF_3\rightarrow 0$$
is an  exact sequence of local systems of finitely generated
right-$R$-modules, $h^{(F_i)_\C}$ are hermitean metrics,  $\eta_i\in
\Omega^{ev}(M)/\image(d)$, and
we form $f_i:=(\cF_i,h^{(F_i)_\C},\eta_i)$, then
$f_2\sim f_1+f_3$ if 
$$\eta_2=\eta_1+\eta_3+\cT(\cC,h^{\cC})\ ,$$
where 
$\cT(\cC,h^{\cC})$ is the analytic torsion form associated to the exact complex
of flat complex vector bundles $(\cC,h^{\cC})$ 
$$0\rightarrow (F_1)_\C\rightarrow (F_2)_\C\rightarrow (F_3)_\C \rightarrow 0$$
equipped with the metric $h^{\cC}$ induced by $h^{(F_i)_\C}$.
For the details of the definition of the torsion form $\cT(\cC,h^\cC)$  we
refer to \cite{bismutlott95}, Sec. 2, or \cite{lott99}, A.3. 
The appearence of the torsion form in the equivalence relation is explained
by the relation
 $$d\cT(\cC,h^{\cC})=\sum_{i=1}^3 (-1)^i c((F_i)_\C,h^{(F_i)_\C})\ . $$

Lott shows that 
$$f \mapsto c(\nabla^{F_\C},h^{F_\C})-d\eta$$
extends to a map $c^\prime:\hat{K}^0_R(M)\rightarrow \Omega^{odd}(M)$,
and he defines
$$\bar{K}^0_R(M):=\ker(c^\prime)\ .$$
The assignment $M\mapsto \bar{K}^0_R(M)$ yields a homotopy invariant
contravariant functor from the category of manifolds to abelian groups.

We now consider a smooth fibre bundle $p:E\rightarrow B$ with compact fibre
$Z$. If $\cF$  is a locally constant sheaf of finitely
generated right-$R$-modules, then we can form the sheaves
$\ccR^ip_*\cF$ on $B$ which are again locally constant sheaves of finitely
generated right-$R$-modules. 

If we choose a fibrewise Riemannian metric $g^Z$ (i.e. a metric on the
vertical bundle $TZ$), then we can compute 
$(\ccR^ip_*\cF)_\C=\ccR^ip_*(\cF_\C)$ (this equality holds because $\C^n$
is a flat $R$-module) using the fibrewise de Rham complex twisted with $F_\C$.
The metric $g^Z$ and a hermitean meric $h^{F_\C}$ induce $L^2$-scalar products.
Identifying the stalk $(\ccR^ip_*(\cF)_\C)_b$ (which is the fibre of the
flat complex vector bundle $(R^ip_*\cF)_\C$) 
with harmonic forms we obtain metrics $h^{(R^ip_*\cF)_\C}$.

We further choose a horizontal distribution $T^HE$. It induces a connection
$\nabla^{TZ}$ on the vertical bundle. Let $e(TZ,\nabla^{TZ})$ be the
associated Euler form. We refer to \cite{bismutlott95} for the definition of
the analytic torsion form $\cT(T^HE,g^Z,h^{F_\C})\in\Omega^{ev}(B)$.

J. Lott defines the push-forward $p_!:\bar{K}^0_R(M)\rightarrow
\bar{K}^0_R(B)$ by the assignment:
$$(\cF,h^F,\eta)\mapsto
\sum_p (\ccR^ip_*(\cF),h^{(R^ip_*\cF)_\C},0) +
(0,0,\int_Ze(TZ,\nabla^{TZ})\wedge \eta-\cT(T^HE,g^Z,h^{F_\C}))\ .$$
Lott proves well-definedness and independence of  $T^HE$ and $g^Z$.

\section{Analytic torsion form and iterated fibrations}\label{ma}

This section reviews the result of \cite{ma99}.
Let $p:E\rightarrow B$ be a fibration with fibre $Z$ which is in fact
an iterated fibration. We assume that there are fibrations 
$p_1:E:=E_1\rightarrow E_2$ with fibre $Z_1$ and
$p_2:E_2\rightarrow B$ with fibre $Z_2$ such that $Z\rightarrow Z_2$ is a
fibration with fibre $Z_1$. By $TZ$, $TZ_1$, and $TZ_2$, we denote the
corresponding vertical bundles. We choose vertical Riemannian metrics $g^Z$,
$g^{Z_1}$, $g^{Z_2}$.

Furthermore, we choose horizontal bundles $T^HE$, $T^HE_1$,
and $T^HE_2$, for $p$, $p_1$, $p_2$.
We obtain connections $\nabla^{TZ}$, $\nabla^{TZ_1}$, and $\nabla^{TZ_2}$.
We identify $TZ$ with $TZ_1\oplus p_1^* TZ_2$ (using $T^HE_1$) and obtain
another connection ${}^0\nabla^{TZ}:=\nabla^{TZ_1}\oplus p_1^*\nabla^{TZ_2}$
on $TZ$. By $\tilde e(TZ,\nabla^{TZ},{}^0\nabla^{TZ})$ we denote the
corresponding transgression of the Euler form such that
$$d\tilde e(TZ,\nabla^{TZ},{}^0\nabla^{TZ})=e(TZ,{}^0\nabla^{TZ})-
e(TZ,\nabla^{TZ})\ .$$

The main result of X. Ma is the formula
\begin{eqnarray*}
&&\cT(T^HE,g^Z,h^{F_\C})-\int_{Z_2}
e(TZ_2,\nabla^{TZ_2})\wedge \cT(T^HE_1,g^{Z_1},h^F)\\
&&-\sum_{i}(-1)^i\cT(T^HE_2,g^{Z_2},h^{(R^ip_*\cF)_\C})\\
&&+\int_{Z}\tilde e(TZ,\nabla^{TZ},{}^0\nabla^{TZ})\wedge 
c(\nabla^{F_\C},h^{F_\C})\\ &=&\cS\:\: (mod\: \image(d))\ ,
\end{eqnarray*}
where $\cS$ is a higher analytic torsion invariant associated to the Leray
spectral sequence of the family of fibrations $Z_b\rightarrow (Z_1)_b$, $b\in
B$.

We now describe $\cS$ in detail. The spectral sequence $(\cE_r,d_r)$ 
is associated to the composition $(p_2)_*\circ (p_1)_*$. Its second term is
$\cE_2^{p,q}=\ccR^p(p_2)_* \ccR^q(p_1)_*(\cF)$, and it converges to
$\ccR^{p+q}p_*(\cF)$. Since $\C^n$ is a flat $R$-module complexification
commutes with taking the spectral sequence. In particular, we obtain
flat complex vector bundles $(E_r^{p,q})_\C$. The differentials $d_r$ of the
spectral sequence induce corresponding bundle homomorphisms
such that we obtain complexes of flat complex vector bundles
$$ \dots\stackrel{d_r}{\rightarrow}
(E_r^{p,q})_\C\stackrel{d_r}{\rightarrow}
(E_r^{p+r,q+1-r})_\C\stackrel{d_r}{\rightarrow}\dots$$ with cohomology
$(E^{p,q}_{r+1})_\C$.  $h^{(R^q(p_*)_1\cF)_\C}$ induces metrics
$h^{(E_2)_\C^{p,q}}$. Now we obtain inductively metrics on the cohomology
groups $h^{(E^{p,q}_{r+1})_\C}$.

In order to save notation we denote by $\cD_r$ the direct sum of complexes
above at the level $r$ and by $h^{\cD_r}$ the induced metric.

Let
$$\cV:\dots\stackrel {d_{i-1}}{\rightarrow} V_i\stackrel{d_i}{\rightarrow}
V_{i+1}\stackrel{d_{i+1}}{\rightarrow}\dots $$ 
be a finite complex of flat complex
vector bundles equipped with hermitean metrics $h^{V_i}$ (we write $h^{\cV}$
for the whole collection). We further choose hermitean metrics on the flat
cohomology bundles $h^{H^i}$ (and we again write $h^{\cH}$ for this
collection).  We form the short exact sequences \begin{eqnarray*}
\cC_i&:& 0\rightarrow \ker(d_i)\rightarrow
V_i\rightarrow \image(d_i)\rightarrow 0\\
\cD_i&:&0\rightarrow \image(d_{i-1}) \rightarrow
\ker(d_i) \rightarrow  H^i \rightarrow 0
\end{eqnarray*}
where all spaces have induced hermitean metrics $h^{\cC_i}$, $h^{\cD_i}$.
We define
$$\cT(\cV,h^{\cV},h^{\cH}):=\sum_{i} (-1)^i \left(
\cT(\cC_i,h^{\cC_i})+\cT(\cD_i,h^{\cD_i})\right)\ .$$

If $V$ is a flat complex vector bundle with a filtration 
$0\subset V_0\subset V_1\subset\dots\subset V_n=V$
by flat subbundles, then
we consider the short exact sequences
$$\cE_i:0\rightarrow V_i \rightarrow V_{i+1}\rightarrow
Gr_{i+1}(V):=V_{i+1}/V_i\rightarrow 0\ .$$
If we further choose hermitean metrics $h^{V}$ (inducing $h^{V_i}$ by
restriction) and $h^{Gr_i(V)}$, then we define metrics $h^{\cE_i}$ and
$$\cT(V,Gr(V),h^V,h^{Gr(V)}):=\sum_i \cT(\cE_i,h^{\cE_i})\ .$$

We can now define
$$\cS:=\sum_{r=2}^\infty
\cT(\cD_r,h^{\cD_r},h^{\cD_{r+1}})-\sum_{k=0}^\infty(-1)^k 
\cT((R^kp_*\cF)_\C,\oplus_{p+q=k}(\cE_\infty^{p,q})_\C,
h^{(R^kp_*\cF)_\C},h^{(\cE_\infty^{p,q})_\C})\ .$$
In the last term we use the natural identification
$Gr_{k-p}(R^kp_*\cF)_\C\stackrel{\cong}{\rightarrow} (\cE_\infty^{p,q})_\C$.

\section{Verification of Theorem \ref{main}}\label{pf}

The group $\bar{K}^0_R(E)$ is generated by elements 
$f:=(\cF,h^{F_\C},\eta)$ with $d\eta=c(\nabla^{F_\C},h^{F_\C})$.
Let us write out a representative of $(p_1)_!([f])$.
We obtain $$\sum_{q}(-1)^q(\ccR^q(p_1)_*\cF,h^{(R^q(p_1)_*\cF)_\C},0)+ 
(0,0,\int_{Z_1}e(TZ_1,\nabla^{TZ_1})\wedge \eta-\cT(T^HE_1,g^{Z_1},h^{F_\C}))\
.$$ A representative of $(p_2)_!\circ (p_1)_!([f])$ is given by
\begin{eqnarray*}
&&\sum_{p,q} (-1)^{p+q}
(\ccR^p(p_2)_*\ccR^q(p_1)_*\cF,h^{(E_2)_\C^{p,q}},0)\\ &&+(0,0,\int_{Z_2}
e(TZ_2,\nabla^{TZ_2})\wedge \left(\int_{Z_1}e(TZ_1,\nabla^{TZ_1})\wedge
\eta-\cT(T^HE_1,g^{Z_1},h^{F_\C})\right))\\
&&-\sum_{q}(-1)^q(0,0,\cT(T^HE_2,g^{Z_2},h^{(R^q(p_1)_*\cF)_\C}))\ .
\end{eqnarray*}
We must show that this expression represents the same element as
$$ 
\sum_i(-1)^i (\ccR^ip_*\cF,h^{(R^ip_*\cF)_\C},0) +
(0,0,\int_Ze(TZ,\nabla^{TZ})\wedge \eta-\cT(T^HE,g^Z,h^{F_\C}))\ .$$
We first compare the terms involving the form $\eta$.
Indeed we have 
\begin{eqnarray*}\lefteqn{\int_{Z_2} e(TZ_2,\nabla^{TZ_2})\wedge \left(
 \int_{Z_1}e(TZ_1,\nabla^{TZ_1})\wedge
\eta \right)- \int_Ze(TZ,\nabla^{TZ})\wedge
\eta}&&\\
 &=& \int_Z d\tilde e(TZ,\nabla^{TZ},{}^0\nabla^{TZ}) \wedge \eta\\
&=& \int_Z \tilde e(TZ,\nabla^{TZ},{}^0\nabla^{TZ}) \wedge d\eta\\
&=& \int_Z \tilde e(TZ,\nabla^{TZ},{}^0\nabla^{TZ}) \wedge
c(\nabla^{F_\C},h^{F_\C})\:(mod\:\image(d))\ .
\end{eqnarray*}
 Thus it remains to show that
\begin{eqnarray*}
 &&\sum_{p,q} (-1)^{p+q}
(\ccR^p(p_2)_*\ccR^q(p_1)_*\cF,h^{(E_2)_\C^{p,q}},0)\\ &&- 
(0,0,\int_{Z_2} e(TZ_2,\nabla^{TZ_2})\wedge  \cT(T^HE_1,g^{Z_1},h^{F_\C}))\\
&&-\sum_{q}(-1)^q(0,0,\cT(T^HE_2,g^{Z_2},h^{(R^q(p_1)_*\cF)_\C}))\\
&&-\sum_p (-1)^p(\ccR^pp_*(\cF),h^{(R^pp_*\cF)_\C},0) \\
&&+(0,0,\cT(T^HE,g^Z,h^{F_\C})\\
&&+ \int_Z \tilde e(TZ,\nabla^{TZ},{}^0\nabla^{TZ}) \wedge
c(\nabla^{F_\C},h^{F_\C})
\end{eqnarray*}
represents the trivial element in $\hat{K}^0_R(B)$.

Let
$$\cV:\dots\stackrel {d_{i-1}}{\rightarrow} \cV_i\stackrel{d_i}{\rightarrow}
\cV_{i+1}\stackrel{d_{i+1}}{\rightarrow} \dots $$ 
be a finite complex of local systems of finitely generated right-$R$-modules
over $B$. We fix hermitean metrics $h^{(V_i)_\C }$ (we write $h^{\cV_\C}$
for the whole collection). 
Since $\C^n$ is a flat $R$-module we can interchange the operation of
complexification and of taking fibrewise cohomology.
We let $H^i$ denote the flat complex vector bundle obtained from the
complexification of the cohomology  sheaves $\cH^i(\cV)$. We further
choose hermitean metrics $h^{H^i}$ (and we 
write $h^{\cH}$ for this collection).  We consider the short exact sequences
\begin{eqnarray*} \cC_i&:& 0\rightarrow \ker(d_i)\rightarrow \cV_i \rightarrow
\image(d_i)\rightarrow 0\\ 
\cD_i&:&0\rightarrow \image(d_{i-1}) \rightarrow
\ker(d_i) \rightarrow  \cH^i(\cV) \rightarrow 0\ ,
\end{eqnarray*}
where the corresponding complexes of flat complex vector bundles  $(\cC_i)_\C$,
$(\cD_i)_\C$ have induced hermitean metrics $h^{(\cC_i)_\C}$, $h^{(\cD_i)_\C}$.
In $\hat{K}^0_R(B)$ we have
\begin{eqnarray*}
\sum_{i}(-1)^i (\cV_i,h^{(V_i)_\C},0)&=&
\sum_{i}(-1)^i
(\cH^i(\cV),h^{H^i},-\cT(C_i,h^{(\cC_i)_\C})-\cT(D_i,h^{(\cD_i)_\C}))\\
&=&  \sum_{i}(-1)^i
(\cH^i(\cV),h^{H^i},0)-(0,0,\cT((\cV)_\C,h^{(\cV)_\C},h^{\cH_\C}))\ .
\end{eqnarray*}
Using this we compute
\begin{eqnarray*}
&&\sum_{p,q} (-1)^{p+q}
(\ccR^p(p_2)_*\ccR^q(p_1)_*\cF,h^{(E_2)_\C^{p,q}},0)\\ &=&\sum_{p,q}
(-1)^{p+q}
(\cE_3^{p,q},h^{(E_3)_\C^{p,q}},0)-(0,0,\cT(\cD_2,h^{\cD_2},h^{\cD_{3}}))\\
&=&\dots\\ &=&\sum_{p,q} (-1)^{p+q}
(\cE_\infty^{p,q},h^{(E_\infty)_\C^{p,q}},0)-\sum_{r=2}^\infty(0,0,\cT(\cD_r,h^{\cD_r},h^{\cD_{r+1}}))\ .
\end{eqnarray*}

Let now  $\cV$ be a local system of finitely generated right-$R$-modules which
is filtered by  local systems of submodules  $0\subset \cV_0\subset
\cV_1\subset\dots\subset \cV_n=\cV$.
We fix a  hermitean metric $h^{V_\C}$ which induces metrics $h^{(V_i)_\C}$.
Furthermore we fix metrics $h^{Gr_i(\cV)_\C}$.
In $\hat{K}^0_R(B)$ we have
$$
(\cV,h^{V_\C},0)=\sum_{i}(Gr_i(\cV),g^{Gr_i(\cV)_\C},0)-(0,0,\cT(V_\C,Gr(\cV)_\C,h^{V_\C},h^{Gr(\cV)_\C}))\ .
$$

Using this observation we further compute 
\begin{eqnarray*}\lefteqn{
\sum_{p+q=k} (-1)^{k}
(\cE_\infty^{p,q},h^{(E_\infty)_\C^{p,q}},0)}&&\\
&=&
\sum_{i}(Gr_i(\ccR^kp_*\cF),h^{Gr_i(R^kp_*\cF)_\C},0)\\
&=&(\ccR^kp_*\cF,h^{(R^kp_*\cF)_\C},0)+(0,0,\cT((R^kp_*\cF)_\C,\oplus_{p+q=k}(\cE_\infty^{p,q})_\C,
h^{(R^kp_*\cF)_\C},h^{Gr(R^kp_*\cF)_\C}))\ .
\end{eqnarray*}

Thus it remains to show that
\begin{eqnarray*}
&&-\sum_{r=2}^\infty \cT(\cD_r,h^{\cD_r},h^{\cD_{r+1}})\\
&&+
\sum_{k=0}^\infty (-1)^k\cT((R^kp_*\cF)_\C,\oplus_{p+q=k}(\cE_\infty^{p,q})_\C,
h^{(R^kp_*\cF)_\C},h^{Gr(R^kp_*\cF)_\C})\\
&&- \int_{Z_2} e(TZ_2,\nabla^{TZ_2})\wedge 
\cT(T^HE_1,g^{Z_1},h^{F_\C})\\
&&-\sum_{p}(-1)^p\cT(T^HE_2,g^{Z_2},h^{(R^p(p_1)_*\cF)_\C})\\
&&+\cT(T^HE,g^Z,h^{F_\C})\\
&& + \int_Z \tilde e(TZ,\nabla^{TZ},{}^0\nabla^{TZ}) \wedge
c(\nabla^{F_\C},h^{F_\C})
\end{eqnarray*}
is an exact form. But this is exactly the assertion of X. Ma.
This finishes the verification. \hB

\bibliographystyle{plain}

\end{document}